\documentclass[
11pt,%
tightenlines,%
twoside,%
onecolumn,%
nofloats,%
nobibnotes,%
nofootinbib,%
superscriptaddress,%
noshowpacs,%
centertags]%
{revtex4}
\usepackage{ljm}
\begin{document}

\titlerunning{Self-similar solutions to model degenerate equations} 
\authorrunning{M.~Ruzhansky,\,A.~Hasanov} 

\title{Self-similar solutions of some model degenerate partial differential equations of the second, third and fourth order}

\author{\firstname{M.}~\surname{Ruzhansky}}
\email[E-mail: ]{michael.ruzhansky@ugent.be,  ruzhansky@gmail.com}
\affiliation{Department of Mathematics, Analysis, Logic and Discrete Mathematics Ghent University, Belgium}
\affiliation{School of Mathematical Sciences, Queen Mary University of London, United Kingdom}

\author{\firstname{A.}~\surname{Hasanov}}
\email[E-mail: ]{anvarhasanov@yahoo.com} \affiliation{Institute of
Mathematics, 81 Mirzo-Ulugbek street, Tashkent, 700170 Uzbekistan}
\firstcollaboration{(Submitted by A.~A.~Editor-name)} 

\received{July 27, 2020;  revised August --, 2020; accepted August
31, 2019}

\thanks{The authors were supported by the FWO Odysseus 1 grant G.0H94.18N: Analysis and Partial Differential Equations. M.R. was also supported by the EPSRC Grant EP/R003025/1 and by the Leverhulme Research Grant RPG-2017-151.}

\begin{abstract} When studying boundary value problems for some partial differential equations arising in applied mathematics, we often have to study the solution of a system of partial differential equations satisfied by hypergeometric functions and find explicit linearly independent solutions for the system. In this study, we construct self-similar solutions of some model degenerate partial differential equations of the second, third, and fourth order. These self-similar solutions are expressed in terms of hypergeometric functions.
\end{abstract}



\maketitle


\section{Introduction and Preliminaries}

A self-similar solution is a solution to some system or equation in which independent variables do not appear independently, but only in combination. To obtain the desired combination of variables for a self-similar solution one often uses the methods of the theory of dimensions. The methods of the theory of dimensions originate from the works of J. Bertrand J., A. Vaschy, subsequently generalised by H. Weyl.

As an example (\cite[p. 113]{Sedov:1}), we consider the problem of the diffusion of vortices in a viscous incompressible fluid under the assumption that the motion of the fluid is plane-parallel and the fluid occupies the entire plane. The motion in question is transient. Suppose that at the initial moment of time   the fluid can potentially move everywhere, except for the pole, which is a trace on the plane of motion of an infinite rectilinear concentrated vortex with circulation. The equation of vortex propagation in this case has the form
\begin{equation}\label{eq 1.1}
\frac{{\partial \Omega }}{{\partial t}} = \nu \left( {\frac{{{\partial ^2}\Omega }}{{\partial {r^2}}} + \frac{1}{r}\frac{{\partial \Omega }}{{\partial r}}} \right),
\end{equation}
where $\Omega $ is the angular velocity of the fluid particles in concentrated circles, and $\nu $ is the coefficient of kinematic viscosity of the fluid. The solution is sought in the form
\begin{equation}\label{eq 1.2}
\Omega \left( {r,t} \right) = \frac{E}{{\nu t}}\psi \left( \xi  \right),\,\,\,\xi  = \frac{{{r^2}}}{{\nu t}}.
\end{equation}
Substituting (\ref{eq 1.2}) into equation (\ref{eq 1.1}), we obtain the solution
\begin{equation}\label{eq 1.3}
\Omega \left( {r,t} \right) = \displaystyle \frac{E}{{\nu t}}A{e^{ -  \frac{{{r^2}}}{{4\nu t}}}},	
\end{equation}
where $A$ is determined from the initial condition of the problem. Thus, (\ref{eq 1.3}) is a self-similar solution to equation (\ref{eq 1.1}).

It is known that to solve applied problems, one needs to set up a mathematical model of the problem under consideration. In many mathematical models, degenerate differential equations appear (especially in gas dynamics, quantum chemistry, in theoretical physics, in the theory of infinitesimal bending of surfaces of revolution, a momentless theory of shells, etc.).

L.D. Landau and E.M. Lifshits in their article \cite{Landau:2} explored the features of the shock wave flow using the Euler-Tricomi equation
$$
{u_{xx}} + {u_{yy}} + \frac{1}{{3y}}{u_y} = 0,
$$
and defined particular solutions of the form
$$
{u_k} = {x^{2k}}F\left( { - k, - k + \frac{1}{2}; - 2k + \frac{5}{6};1 + \frac{{4{y^2}}}{{9{x^2}}}} \right),
$$ where $k =  \pm \displaystyle\frac{n}{2},\,\,\frac{1}{3} \pm \frac{n}{2},\,\,n \in N_0 $.     

Note that the energy absorbed by a non-ferromagnetic conducting sphere placed in an external inhomogeneous magnetic field is calculated explicitly using the hypergeometric functions of many variables (\cite{Lohofer:3}).

The hypergeometric functions of Kampe de Feriet also appear in theoretical physics and quantum chemistry (see e.g. \cite{Niukkanen:4}). In the monographs \cite{Courant:5} - \cite{Frankl:7}, attention was drawn to the fact that many problems of supersonic gas dynamics are solved using hypergeometric functions.

Using the method of self-similar solutions in articles \cite{Ergashev:8} - \cite{Srivastava:20}, fundamental solutions were found and in articles \cite{Salakhitdinov:21} - \cite{Salakhitdinov:24} the main boundary-value problems for the generalised axisymmetric Helmholtz equation were solved.

In this paper, using the method of self-similar solutions, we construct some special solutions of degenerate partial differential equations that are expressed by hypergeometric functions.

\section{A parabolic equation with one line of degeneration}

Consider in the domain $\Omega  = \left\{ {\left( {x,t} \right):\,\,x > 0,\,\,t > 0} \right\}$, the degenerate parabolic equation
\begin{equation}\label{eq 2.1}
Lu \equiv {u_t} - {u_{xx}} - \frac{{2\alpha }}{x}{u_x} = 0,\,\,\alpha  = const>0.
\end{equation}
We seek self-similar solutions of equation (\ref{eq 2.1}) in the form
\begin{equation}\label{eq 2.2}
u = P\omega \left( \sigma  \right),
\end{equation}
where $\omega=\omega \left( \sigma  \right)$   is an unknown function, and where $\sigma  =  -  \frac{{{x^2}}}{{4t}},\,\,P = {t^{ -  \frac{1}{2}}}$   . Substituting (\ref{eq 2.2}) into equation (\ref{eq 2.1}), we have
\begin{equation}\label{eq 2.3}
P{\omega _{\sigma \sigma }}\sigma _x^2 + \left[ {2{P_x}{\sigma _x} + P\left( {{\sigma _{xx}} + \frac{{2\alpha }}{x}{\sigma _x} - {\sigma _t}} \right)} \right]{\omega _\sigma } + \left( {{P_{xx}} + \frac{{2\alpha }}{x}{P_x} - {P_t}} \right)\omega  = 0.
\end{equation}
After elementary calculations, we find
$$
\sigma _x^2 =  - \frac{1}{t}\sigma ,\,\,2{P_x}{\sigma _x} + P\left( {{\sigma _{xx}} + \frac{{2\alpha }}{x}{\sigma _x} - {\sigma _t}} \right) =  - P\frac{1}{t}\left( {\frac{{1 + 2\alpha }}{2} - \sigma } \right),
$$
$$
{P_{xx}} + \frac{{2\alpha }}{x}{P_x} - {P_t} = \frac{1}{2}P\frac{1}{t}.
$$
Therefore, in view of the indicated equalities, the ordinary differential equation (\ref{eq 2.3}) has the form
\begin{equation}\label{eq 2.4}
\sigma {\omega _{\sigma \sigma }} + \left( {\frac{{1 + 2\alpha }}{2} - \sigma } \right){\omega _\sigma } - \frac{1}{2}\omega  = 0.                  \end{equation}
It is known (\cite{Erdelyi:25}) that the equation
\begin{equation}\label{eq 2.5}
x{w_{xx}} + \left( {c - x} \right){w_x} - aw = 0,
\end{equation}
has two linearly independent solutions
\begin{equation}\label{eq 2.6}
\begin{array}{l}
{w_1} = {c_1}{}_1{F_1}\left( {a;c;x} \right) = {c_1}{e^x}{}_1{F_1}\left( {c - a;c; - x} \right),\\
{w_2} = {c_2}{x^{1 - c}}{}_1{F_1}\left( {a - c + 1;2 - c;x} \right) = {c_2}{x^{1 - c}}{e^x}{}_1{F_1}\left( {1 - a;2 - c; - x} \right),
\end{array}
\end{equation}
where the  hypergeometric function ${}_1{F_1}\left( {a;c;x} \right)$ has the form
$$
{}_1{F_1}\left( {a;c;x} \right) = \sum\limits_{m = 0}^\infty  {} \frac{{{{\left( a \right)}_m}}}{{{{\left( c \right)}_m}m!}}{x^m},
$$
and $${\left( a \right)_m} = \Gamma \left( {a + m} \right)/\Gamma \left( a \right) = a\left( {a + 1} \right)\left( {a + 2} \right) \cdot  \cdot  \cdot \left( {a + m - 1} \right)$$   is the Pochhammer symbol (\cite{Erdelyi:25}).
Therefore, taking into account (\ref{eq 2.6}) and (\ref{eq 2.4}), we define
\begin{equation}\label{eq 2.7}
{\omega _1} = {c_1}{}_1{F_1}\left( {\frac{1}{2};\frac{{1 + 2\alpha }}{2};\sigma } \right),\,\,\,{\omega _2} = {c_2}{\sigma ^{\frac{{1 - 2\alpha }}{2}}}{}_1{F_1}\left( {1 - \alpha ;\frac{{3 - 2\alpha }}{2};\sigma } \right).
\end{equation}
Substituting (\ref{eq 2.7}) into (\ref{eq 2.2}), we finally obtain
\begin{equation}\label{eq 2.8}
{u_1}\left( {x,t} \right) = {c_1}\frac{1}{{\sqrt t }}{}_1{F_1}\left( {\frac{1}{2};\frac{{1 + 2\alpha }}{2}; - \frac{{{x^2}}}{{4t}}} \right),
\end{equation}
\begin{equation}\label{eq 2.9}
{u_2}\left( {x,t} \right) = {c_2}\displaystyle \frac{1}{{\sqrt t }}{\left( {\frac{{{x^2}}}{{4t}}} \right)^{ \frac{{1 - 2\alpha }}{2}}}{}_1{F_1}\left( {1 - \alpha ;\frac{{3 - 2\alpha }}{2}; - \frac{{{x^2}}}{{4t}}} \right),
\end{equation}
two self-similar solutions of equation (\ref{eq 2.1}), where ${c_1},{c_2} $ are constants. Note that in \cite{Karimov:26}, the boundary value problems for the equation ${L^m}u = 0$    were considered.

\section{A parabolic equation with two lines of degeneracy}

In the domain $\Omega  = \left\{ {\left( {x,y,t} \right):\,\,x > 0,\,y > 0,\,\,t > 0} \right\}$, we consider the equation
\begin{equation}\label{eq 3.1}
Lu \equiv {u_t} - {u_{xx}} - {u_{yy}} - \frac{{2\alpha }}{x}{u_x} - \frac{{2\beta }}{y}{u_y} = 0,\,\,\alpha ,\beta  = const.
\end{equation}
The solution to equation (\ref{eq 3.1}) is sought in the form
\begin{equation}\label{eq 3.2}
u = P\omega \left( {\xi ,\eta } \right),                                                                                                                                 \end{equation}
where  $\xi  =  - \displaystyle \frac{{{x^2}}}{{8t}},$ $\eta  =  -\displaystyle \frac{{{y^2}}}{{8t}},$  $P = {t^{ - \frac{1}{2}}},$ and $\omega \left( {\xi ,\eta } \right)$ is an unknown function. Substituting (\ref{eq 3.2}) into (\ref{eq 3.1}), we have
\begin{equation}\label{eq 3.3}
{A_1}{\omega _{\xi \xi }} + {A_2}{\omega _{\xi \eta }} + {A_3}{\omega _{\eta \eta }} + {A_4}{\omega _\xi } + {A_5}{\omega _\eta } + {A_6}\omega  = 0,                                                                                                                                        \end{equation}
where
$$
\begin{array}{l}
{A_1} =\displaystyle P\left( {\xi _x^2 + \xi _y^2} \right),\,\,\,{A_2} = 2P\left( {{\xi _x}{\eta _x} + {\xi _y}{\eta _y}} \right),\,\,\,
{A_3} = P\left( {\eta _x^2 + \eta _y^2} \right),\\
{A_4} = \displaystyle \frac{{2\alpha }}{x}P{\xi _x} + \frac{{2\beta }}{y}P{\xi _y} - P{\xi _t} + 2{P_x}{\xi _x} + P{\xi _{xx}} + 2{P_y}{\xi _y} + P{\xi _{yy}},\\
{A_5} =\displaystyle \frac{{2\alpha }}{x}P{\eta _x} + \frac{{2\beta }}{y}P{\eta _y} - P{\eta _t} + 2{P_x}{\eta _x} + P{\eta _{xx}} + 2{P_y}{\eta _y} + P{\eta _{yy}},\\
{A_6} =  - {P_t} + {P_{xx}} + {P_{yy}} + \displaystyle \frac{{2\alpha }}{x}{P_x} + \displaystyle \frac{{2\beta }}{y}{P_y}.
\end{array}
$$
Calculating the values of the coefficients of equation (\ref{eq 3.3}), we obtain a system of partial differential equations
\begin{equation}\label{eq 3.4}
\left\{ {\begin{array}{*{20}{c}}
{\xi {\omega _{\xi \xi }} + \left( {\displaystyle \frac{{1 + 2\alpha }}{2} - \xi } \right){\omega _\xi } - \eta {\omega _\eta } - \displaystyle \frac{1}{2}\omega  = 0}\\
{\eta {\omega _{\eta \eta }} + \left( {\displaystyle \frac{{1 + 2\beta }}{2} - \eta } \right){\omega _\eta } - \xi {\omega _\xi } - \displaystyle \frac{1}{2}\omega  = 0.}
\end{array}} \right.	
\end{equation}
In the monograph \cite{Erdelyi:25}, \cite{Appell:27} the following system of hypergeometric equations was considered
\begin{equation}\label{eq 3.5}
\left\{ {\begin{array}{*{20}{c}}
{x{w_{xx}} + \left( {{c_1} - x} \right){w_x} - y{w_y} - aw = 0}\\
{y{w_{yy}} + \left( {{c_2} - y} \right){w_y} - x{w_x} - aw = 0}
\end{array}} \right.
\end{equation}
and 4 linearly independent solutions were found, expressed in terms of the confluent Kummer functions,
\begin{equation}\label{eq 3.6}
{w_1} = {\lambda _1}{\Psi _2}\left( {a;{c_1},{c_2};x,y} \right),
\end{equation}
\begin{equation}\label{eq 3.7}
 {w_2} = {\lambda _2}{x^{1 - {c_1}}}{\Psi _2}\left( {a + 1 - {c_1};2 - {c_1},{c_2};x,y} \right),
\end{equation}
\begin{equation}\label{eq 3.8}
{w_3} = {\lambda _3}{y^{1 - {c_2}}}{\Psi _2}\left( {a + 1 - {c_2};{c_1},2 - {c_2};x,y} \right),
\end{equation}
\begin{equation}\label{eq 3.9}
{w_4} = {\lambda _4}{x^{1 - {c_1}}}{y^{1 - {c_2}}}{\Psi _2}\left( {a + 2 - {c_1} - {c_2};2 - {c_1},2 - {c_2};x,y} \right),
\end{equation}
where ${\lambda _i} = const,\,\,i = {1,2,3,4},$   and
$$
{\Psi _2}\left( {a;{c_1},{c_2};x,y} \right) = \sum\limits_{m,n = 0}^\infty  {} \frac{{{{\left( a \right)}_{m + n}}}}{{{{\left( {{c_1}} \right)}_m}{{\left( {{c_2}} \right)}_n}m!n!}}{x^m}{y^n}.
$$
In view of (\ref{eq 3.6}) - (\ref{eq 3.9}), we define
$$
{\omega _1} = {\lambda _1}{\Psi _2}\left( {\displaystyle \frac{1}{2};\frac{{1 + 2\alpha }}{2},\frac{{1 + 2\beta }}{2}; - \frac{{{x^2}}}{{8t}}, - \frac{{{y^2}}}{{8t}}} \right),
$$
$$
{\omega _2} = {\lambda _2}{\left( {\displaystyle \frac{{{x^2}}}{{8t}}} \right)^{ \frac{{1 - 2\alpha }}{2}}}{\Psi _2}\left( {1 - \alpha ;\frac{{3 - 2\alpha }}{2},\frac{{1 + 2\beta }}{2}; - \frac{{{x^2}}}{{8t}}, - \frac{{{y^2}}}{{8t}}} \right),
$$
$$
{\omega _3} = {\lambda _3}{\left( {\frac{{{y^2}}}{{8t}}} \right)^{ \frac{{1 - 2\beta }}{2}}}{\Psi _2}\left( {1 - \beta ;\frac{{1 + 2\alpha }}{2},\frac{{3 - 2\beta }}{2}; - \frac{{{x^2}}}{{8t}}, - \frac{{{y^2}}}{{8t}}} \right),
$$
$$
{\omega _4} = {\lambda _4}{\left( {\frac{{{x^2}}}{{8t}}} \right)^{ \frac{{1 - 2\alpha }}{2}}}{\left( {\frac{{{y^2}}}{{8t}}} \right)^{ \frac{{1 - 2\beta }}{2}}}{\Psi _2}\left( {\frac{{3 - 2\alpha  - 2\beta }}{2};\frac{{3 - 2\alpha }}{2},\frac{{3 - 2\beta }}{2}; - \frac{{{x^2}}}{{8t}}, - \frac{{{y^2}}}{{8t}}} \right).
$$
Substituting  ${\omega _i},\,\,i =  {1,2,3,4}, $  in (\ref{eq 3.2}), we obtain the following special solutions of equation (\ref{eq 3.1}):
\begin{equation}\label{eq 3.10}
{u_1}\left( {x,y,t} \right) = {\lambda _1}\displaystyle \frac{1}{{\sqrt t }}{\Psi _2}\left( {\frac{1}{2};\frac{{1 + 2\alpha }}{2},\frac{{1 + 2\beta }}{2}; - \frac{{{x^2}}}{{8t}}, - \frac{{{y^2}}}{{8t}}} \right),
\end{equation}
\begin{equation}\label{eq 3.11}
{u_2}\left( {x,y,t} \right) = {\lambda _2}\displaystyle \frac{{{\displaystyle x^{1 - 2\alpha }}}}{{{t^{1 - \alpha }}}}{\Psi _2}\left( {1 - \alpha ;\frac{{3 - 2\alpha }}{2},\frac{{1 + 2\beta }}{2}; - \frac{{{x^2}}}{{8t}}, - \frac{{{y^2}}}{{8t}}} \right),
\end{equation}
\begin{equation}\label{eq 3.12}
{u_3}\left( {x,y,t} \right) = {\lambda _3}\frac{{{y^{1 - 2\beta }}}}{{{t^{1 - \beta }}}}{\Psi _2}\left( {1 - \beta ;\frac{{1 + 2\alpha }}{2},\frac{{3 - 2\beta }}{2}; - \frac{{{x^2}}}{{8t}}, - \frac{{{y^2}}}{{8t}}} \right),
\end{equation}
\begin{equation}\label{eq 3.13}
{u_4}\left( {x,y,t} \right) = {\lambda _4}\frac{{{x^{1 - 2\alpha }}{y^{1 - 2\beta }}}}{{{t^{2 - \alpha  - \beta }}}}{\Psi _2}\left( {\frac{{3 - 2\alpha  - 2\beta }}{2};\frac{{3 - 2\alpha }}{2},\frac{{3 - 2\beta }}{2}; - \frac{{{x^2}}}{{8t}}, - \frac{{{y^2}}}{{8t}}} \right),
\end{equation}
where ${\lambda _1},{\lambda _2},{\lambda _3},{\lambda _4}$ are constants.

\section{A differential equation of the third order with one line of degeneration}

Consider the equation
\begin{equation}\label{eq 4.1}
Lu \equiv {y^m}{u_{xxx}} - {u_{yyy}} = 0,\,\,m = const > 0,
\end{equation}
in the domain of $\Omega  = \left\{ {\left( {x,y} \right):\,\, - \infty  < x <  + \infty ,y > 0} \right\}$. Special solutions of equation (\ref{eq 4.1}) are sought in the form
\begin{equation}\label{eq 4.2}
u = P\omega \left( \sigma  \right),
\end{equation}
where
$$
P = {x^{ - 3}},\,\,\sigma  = {\left( { - \displaystyle \frac{3}{{x\left( {m + 3} \right)}}{y^{ \frac{{m + 3}}{3}}}} \right)^3},\,\,\beta  = \frac{m}{{m + 3}}.
$$
Substituting (\ref{eq 4.2}) into equation (\ref{eq 4.1}), we find
\begin{equation}\label{eq 4.3}
A{\omega _{\sigma \sigma \sigma }} + B{\omega _{\sigma \sigma }} + C{\omega _\sigma } + D\omega  = 0,
\end{equation}
where
$$A = P\left( {{y^m}\sigma _x^3 - \sigma _y^3} \right),$$
$$B = 3\left[ {{y^m}{P_x}\sigma _x^2 - {P_y}\sigma _y^2 + P\left( {{y^m}{\sigma _x}{\sigma _{xx}} - {\sigma _y}{\sigma _{yy}}} \right)} \right],$$
$$C = \left[ {3\left( {{y^m}{P_{xx}}{\sigma _x} - {P_{yy}}{\sigma _y}} \right) + 3\left( {{y^m}{P_x}{\sigma _{xx}} - {P_y}{\sigma _{yy}}} \right) + P\left( {{y^m}{\sigma _{xxx}} - {\sigma _{yyy}}} \right)} \right],$$
$$D = {y^m}{P_{xxx}} - {P_{yyy}}.$$
After elementary calculations, we have
$$A = \displaystyle \frac{{{3^3}P{y^m}}}{{{x^3}}}{\sigma ^2}\left( {1 - \sigma } \right),$$
$$B = \frac{{{3^3}{y^m}P}}{{{x^3}}}\left[ {\frac{{2 + \beta }}{3} + \frac{{1 + 2\beta }}{3} + 1 - \left( {3 + 1 + \frac{4}{3} + \frac{5}{3}} \right)\sigma } \right]\sigma,$$
$$C = \frac{{{3^3}{y^m}P}}{{{x^3}}}\left[ {\frac{{2 + \beta }}{3}\frac{{1 + 2\beta }}{3} - \left( {1 + 1 + \frac{4}{3} + \frac{5}{3} + 1 \cdot \frac{4}{3} + 1 \cdot \frac{5}{3} + \frac{4}{3} \cdot \frac{5}{3}} \right)\sigma } \right],$$
$$D =  - \frac{{{3^3}{y^m}P}}{{{x^3}}} \cdot 1 \cdot \frac{4}{3} \cdot \frac{5}{3}.$$
By the above equalities, from (\ref{eq 4.3}) it follows that
\begin{equation}\label{eq 4.4}
\begin{array}{l}
{\sigma ^2}\left( {1 - \sigma } \right){\omega _{\sigma \sigma \sigma }} + \left[ {\displaystyle \frac{{2 + \beta }}{3} + \displaystyle \frac{{1 + 2\beta }}{3} + 1 - \left( {3 + 1 + \frac{4}{3} + \frac{5}{3}} \right)\sigma } \right]\sigma {\omega _{\sigma \sigma }}  \\
 + \left[ {\displaystyle \frac{{2 + \beta }}{3}\frac{{1 + 2\beta }}{3} - \left( {1 + 1 + \frac{4}{3} + \frac{5}{3} + 1 \cdot \frac{4}{3} + 1 \cdot \frac{5}{3} + \frac{4}{3} \cdot \frac{5}{3}} \right)\sigma } \right]{\omega _\sigma } - 1 \cdot \displaystyle \frac{4}{3} \cdot \frac{5}{3}\omega  = 0.
\end{array}
\end{equation}
Thus, we have obtained the ordinary Clausen differential equation (\cite{Appell:27}), which has the form
\begin{equation}\label{eq 4.5}
\begin{array}{l}
{x^2}\left( {1 - x} \right){w_{xxx}} + \left[ {{c_1} + {c_2} + 1 - \left( {3 + {a_1} + {a_2} + {a_3}} \right)x} \right]x{w_{xx}} \\
 + \left[ {{c_1}{c_2} - \left( {1 + {a_1} + {a_2} + {a_3} + {a_1}{a_2} + {a_1}{a_3} + {a_2}{a_3}} \right)x} \right]{w_x} - {a_1}{a_2}{a_3}w = 0.
\end{array}
\end{equation}
The Clausen equation (\ref{eq 4.5}) has the following three linearly independent solutions \cite{Appell:27}:
\begin{equation}\label{eq 4.6}
{w_1} = {\lambda _1}{}_3{F_2}\left[ {\begin{array}{*{20}{c}}
{{a_1},{a_2},{a_3}}\\
{{c_1},{c_2}}
\end{array}x} \right],
\end{equation}
\begin{equation}\label{eq 4.7}
{w_2} = {\lambda _2}{x^{1 - {c_1}}}{}_3{F_2}\left[ {\begin{array}{*{20}{c}}
{{a_1} + 1 - {c_1},{a_2} + 1 - {c_1},{a_3} + 1 - {c_1}}\\
{2 - {c_1},{c_2} + 1 - {c_1}}
\end{array}x} \right],
\end{equation}
\begin{equation}\label{eq 4.8}
{w_3} = {\lambda _3}{x^{1 - {c_2}}}{}_3{F_2}\left[ {\begin{array}{*{20}{c}}
{{a_1} + 1 - {c_2},{a_2} + 1 - {c_2},{a_3} + 1 - {c_2}}\\
{{c_1} + 1 - {c_2},2 - {c_2}}
\end{array}x} \right].
\end{equation}
Note that the Clausen function can be represented as
$$
\begin{array}{l}
{}_3{F_2}\left[ {\begin{array}{*{20}{c}}
{{a_1},{a_2},{a_3}}\\
{{c_1},{c_2}}
\end{array}x} \right] = \displaystyle \frac{{\Gamma \left( {{c_1}} \right)\Gamma \left( {{c_2}} \right)}}{{\Gamma \left( {{a_1}} \right)\Gamma \left( {{a_2}} \right)\Gamma \left( {{c_1} - {a_1}} \right)\Gamma \left( {{c_2} - {a_2}} \right)}}\\
 \times \displaystyle \int\limits_0^1 {\int\limits_0^1 {} {\xi ^{{a_1} - 1}}{\eta ^{{a_2} - 1}}{{\left( {1 - \xi } \right)}^{{c_1} - {a_1} - 1}}} {\left( {1 - \eta } \right)^{{c_2} - {a_2} - 1}}{\left( {1 - x\xi \eta } \right)^{ - {a_3}}}d\xi d\eta ,\\
{\mathop{\rm Re}\nolimits} \,\,{c_i}\,\, > \,{\mathop{\rm Re}\nolimits} \,\,{a_i}\, > 0,\,\,i = 1,2.,
\end{array}
$$
$$
{}_3{F_2}\left( {{a_1},{a_2},{a_3};{c_1},{c_2};x} \right) = {}_3{F_2}\left[ {\begin{array}{*{20}{c}}
{{a_1},{a_2},{a_3}}\\
{{c_1},{c_2}}
\end{array}x} \right] =\displaystyle \sum\limits_{m = 0}^\infty  {} \displaystyle \frac{{{{\left( {{a_1}} \right)}_m}{{\left( {{a_2}} \right)}_m}{{\left( {{a_3}} \right)}_m}}}{{{{\left( {{c_1}} \right)}_m}{{\left( {{c_2}} \right)}_m}m!}}{x^m}.
$$
Taking into account (\ref{eq 4.6}) - (\ref{eq 4.8}), it follows from (\ref{eq 4.4}) that
\begin{equation}\label{eq 4.9}
{\omega _1}\left( \sigma  \right) = {\lambda _1}{}_3{F_2}\left[ {\begin{array}{*{20}{c}}\displaystyle
{1,\frac{4}{3},\frac{5}{3};}\\
{\displaystyle \frac{{2 + \beta }}{3},\frac{{1 + 2\beta }}{3};}
\end{array}\sigma } \right],
\end{equation}
\begin{equation}\label{eq 4.10}
{\omega _2}\left( \sigma  \right) = {\lambda _2}{\sigma ^{ \frac{{1 - \beta }}{3}}}{}_3{F_2}\left[ {\begin{array}{*{20}{c}}
\displaystyle {\displaystyle \frac{{4 - \beta }}{3},\displaystyle \frac{{5 - \beta }}{3},\frac{{6 - \beta }}{3};}\\
\displaystyle {\frac{{4 - \beta }}{3},\frac{{2 + \beta }}{3};}
\end{array}\sigma } \right] = {\lambda _2}{\sigma ^{\frac{{1 - \beta }}{3}}}F\left( {\frac{{5 - \beta }}{3},\frac{{6 - \beta }}{3};\frac{{2 + \beta }}{3};\sigma } \right),
\end{equation}
\begin{equation}\label{eq 4.11}
{\omega _3}\left( \sigma  \right) = {\lambda _3}{\sigma ^{ \frac{{2 - 2\beta }}{3}}}{}_3{F_2}\left[ {\begin{array}{*{20}{c}}
{\displaystyle \frac{{5 - 2\beta }}{3},\displaystyle \frac{{6 - 2\beta }}{3},\frac{{7 - 2\beta }}{3};}\\
{\displaystyle \frac{{4 - \beta }}{3},\displaystyle \frac{{5 - 2\beta }}{3};}
\end{array}\sigma } \right] = {\lambda _3}{\sigma ^{ \frac{{2 - 2\beta }}{3}}}F\left( {\displaystyle \frac{{6 - 2\beta }}{3},\displaystyle \frac{{7 - 2\beta }}{3};\frac{{4 - \beta }}{3};\sigma } \right),
\end{equation}
where $F\left( {a,b;c;x} \right)$    is the Gauss hypergeometric function (\cite{Erdelyi:25}, \cite{Appell:27}). Substituting (\ref{eq 4.9})  - (\ref{eq 4.11}), we finally find the following special solutions to equation (\ref{eq 4.1}):
\begin{equation}\label{eq 4.12}
{u_1}\left( {x,y} \right) = {\lambda _1}{x^{ - 3}}{}_3{F_2}\left[ {\begin{array}{*{20}{c}}\displaystyle
{1,\frac{4}{3},\frac{5}{3};}\\
\displaystyle {\frac{{2 + \beta }}{3},\frac{{1 + 2\beta }}{3};}
\end{array}{{\left( { - \frac{3}{{x\left( {m + 3} \right)}}{y^{ \frac{{m + 3}}{3}}}} \right)}^3}} \right],
\end{equation}
\begin{equation}\label{eq 4.13}
{u_2}\left( {x,y} \right) = {\lambda _2}{x^{\beta  - 4}}yF\left( {\frac{{5 - \beta }}{3},\frac{{6 - \beta }}{3};\frac{{2 + \beta }}{3};{{\left( { - \frac{3}{{x\left( {m + 3} \right)}}{y^{ \frac{{m + 3}}{3}}}} \right)}^3}} \right),
\end{equation}
\begin{equation}\label{eq 4.14}
{u_3}\left( {x,y} \right) = {\lambda _3}{x^{2\beta  - 5}}{y^2}F\left( {\frac{{6 - 2\beta }}{3},\frac{{7 - 2\beta }}{3};\frac{{4 - \beta }}{3};{{\left( { - \frac{3}{{x\left( {m + 3} \right)}}{y^{ \frac{{m + 3}}{3}}}} \right)}^3}} \right),
\end{equation}
where ${\lambda _1},{\lambda _2},{\lambda _3}$  are arbitrary constants.

\section{The third-order differential equation of three variables}

In the domain $\Omega  = \left\{ {\left( {x,y,t} \right):\,\,x > 0,\,y > 0,\,t > 0} \right\}$, we consider the equation
\begin{equation}\label{eq 5.1}
Lu \equiv {x^n}{y^m}\,{u_t} - {t^k}{y^m}{u_{xxx}} - {t^k}{x^n}{u_{yyy}} = 0,\,\,m,n,k = const > 0.
\end{equation}
The solution to equation (\ref{eq 5.1}) is sought in the form
\begin{equation}\label{eq 5.2}
u = P\omega \left( {\xi ,\eta } \right),
\end{equation}
where
$$
P = {\left( {\frac{2}{{k + 1}}{t^{k + 1}}} \right)^{ - 1}},\,\,\xi  =  - \frac{{k + 1}}{{2{{\left( {n + 3} \right)}^3}{t^{k + 1}}}}{x^{n + 3}},\,\,\eta  =  - \frac{{k + 1}}{{2{{\left( {m + 3} \right)}^3}{t^{k + 1}}}}{y^{m + 3}},
$$
$$\alpha  = n/\left( {n + 3} \right),\,\,\beta  = m/\left( {m + 3} \right).$$
Then substituting (\ref{eq 5.2}) into (\ref{eq 5.1}), we obtain a third-order partial differential equation
\begin{equation}\label{eq 5.3}
{A_1}{\omega _{\xi \xi \xi }} + {A_2}{\omega _{\eta \eta \eta }} + {A_3}{\omega _{\xi \xi \eta }} + {A_4}{\omega _{\xi \eta \eta }} + {A_5}{\omega _{\xi \xi }} + {A_6}{\omega _{\xi \eta }} + {A_7}{\omega _{\eta \eta }} + {A_8}{\omega _\xi } + {A_9}{\omega _\eta } + {A_{10}}\omega  = 0,
\end{equation}
where
$$
{A_1} = P{t^k}\left( {{y^m}\xi _x^3 + {x^n}\xi _y^3} \right),\,\,{A_2} = P{t^k}\left( {{y^m}\eta _x^3 + {x^n}\eta _y^3} \right),\,\,{A_3} = 3{t^k}P\left( {{y^m}\xi _x^2{\eta _x} + {x^n}\xi _y^2{\eta _y}} \right),
$$
$$
{A_4} = 3{t^k}P\left( {{y^m}{\xi _x}\eta _x^2 + {x^n}{\xi _y}\eta _y^2} \right),\,\,{A_5} = 3{t^k}\left[ {{y^m}{P_x}\xi _x^2 + {x^n}{P_y}\xi _y^2 + P\left( {{y^m}{\xi _x}{\xi _{xx}} + {x^n}{\xi _y}{\xi _{yy}}} \right)} \right],
$$
$$
{A_6} = 3{t^k}\left[ {2\left( {{y^m}{P_x}{\xi _x}{\eta _x} + {x^n}{P_y}{\xi _y}{\eta _y}} \right) + P\left( {{y^m}{\xi _{xx}}{\eta _x} + {x^n}{\xi _{yy}}{\eta _y}} \right) + P\left( {{y^m}{\xi _x}{\eta _{xx}} + {x^n}{\xi _y}{\eta _{yy}}} \right)} \right],
$$
$$
{A_7} = 3{t^k}\left[ {{y^m}{P_x}\eta _x^2 + {x^n}{P_y}\eta _y^2 + P\left( {{y^m}{\eta _x}{\eta _{xx}} + {x^n}{\eta _y}{\eta _{yy}}} \right)} \right],
$$
$$
\begin{array}{l}
{A_8} = 3{t^k}\left( {{y^m}{P_{xx}}{\xi _x} + {x^n}{P_{yy}}{\xi _y}} \right) + 2{t^k}\left( {{y^m}{P_x}{\xi _{xx}} + {x^n}{P_y}{\xi _{yy}}} \right) + {t^k}\left( {{y^m}{P_x}{\xi _{xx}} + {x^n}{P_y}{\xi _{yy}}} \right)\\
\,\,\,\,\,\,\,\,\, + {t^k}P\left( {{y^m}{\xi _{xxx}} + {x^n}{\xi _{yyy}}} \right) - {x^n}{y^m}P{\xi _t},
\end{array}
$$
$$
\begin{array}{l}
{A_9} = 3{t^k}\left( {{y^m}{P_{xx}}{\eta _x} + {x^n}{P_{yy}}{\eta _y}} \right) + 2{t^k}\left( {{y^m}{P_x}{\eta _{xx}} + {x^n}{P_y}{\eta _{yy}}} \right) + {t^k}\left( {{y^m}{P_x}{\eta _{xx}} + {x^n}{P_y}{\eta _{yy}}} \right)\\
\,\,\,\,\,\,\, \,\,\,+ {t^k}P\left( {{y^m}{\eta _{xxx}} + {x^n}{\eta _{yyy}}} \right) - {x^n}{y^m}P{\eta _t},
\end{array}
$$
$${A_{10}} = {t^k}{y^m}{P_{xxx}} + {t^k}{x^n}{P_{yyy}} - {x^n}{y^m}{P_t}.$$
After some calculations, we have
$${A_1} =  - P{t^k}{x^n}{y^m}\frac{{k + 1}}{{2{t^{k + 1}}}}{\xi ^2},\,\,{A_2} =  - P{t^k}{x^n}{y^m}\frac{{k + 1}}{{2{t^{k + 1}}}}{\eta ^2},\,\,{A_3} = 0,\,\,{A_4} = 0,$$
$${A_5} =  - P{t^k}{x^n}{y^m}\left( {\alpha  + 2} \right)\frac{{k + 1}}{{2{t^{k + 1}}}}\xi ,\,\,{A_6} = 0,\,\,{A_7} =  - P{t^k}{x^n}{y^m}\left( {2 + \beta } \right)\frac{{k + 1}}{{2{t^{k + 1}}}}\eta ,$$
$${A_8} =  - P{t^k}{x^n}{y^m}\frac{{k + 1}}{{2{t^{k + 1}}}}\left( {\frac{{2 + \alpha }}{3}\frac{{1 + 2\alpha }}{3} - 2\xi } \right),\,\,{A_9} =  - P{t^k}{x^n}{y^m}\frac{{k + 1}}{{2{t^{k + 1}}}}\left( {\frac{{2 + \beta }}{3}\frac{{1 + 2\beta }}{3} - 2\eta } \right),$$
$${A_{10}} = P{t^k}{x^n}{y^m}\frac{{k + 1}}{{{t^{k + 1}}}}.$$
Therefore, using the indicated equalities from (\ref{eq 5.3}), we define
\begin{equation}\label{eq 5.4}
\left\{ {\begin{array}{*{20}{c}}
{{\xi ^2}{\omega _{\xi \xi \xi }} + \left( {\displaystyle \frac{{2 + \alpha }}{3} + \displaystyle \frac{{1 + 2\alpha }}{3} + 1} \right)\xi {\omega _{\xi \xi }} + \left( \displaystyle {\frac{{2 + \alpha }}{3}\frac{{1 + 2\alpha }}{3} - \xi } \right){\omega _\xi } - \eta {\omega _\eta } - \omega  = 0}\\
{{\eta ^2}{\omega _{\eta \eta \eta }} + \left( \displaystyle {\frac{{2 + \beta }}{3} + \frac{{1 + 2\beta }}{3} + 1} \right)\eta {\omega _{\eta \eta }} + \left( \displaystyle {\frac{{2 + \beta }}{3}\frac{{1 + 2\beta }}{3} - \eta } \right){\omega _\eta } - \xi {\omega _\xi } - \omega  = 0.}
\end{array}} \right.	
\end{equation}
From the general theory it is easy to determine that the system of hypergeometric equations
\begin{equation}\label{eq 5.5}
\left\{ {\begin{array}{*{20}{c}}
{{x^2}{w_{xxx}} + \left( {{c_2} + {c_1} + 1} \right)x{w_{xx}} + \left( {{c_1}{c_2} - x} \right){w_x} - y{w_y} - aw = 0}\\
{{y^2}{w_{yyy}} + \left( {{d_2} + {d_1} + 1} \right)y{w_{yy}} + \left( {{d_1}{d_2} - y} \right){w_y} - x{w_x} - aw = 0}
\end{array}} \right.
\end{equation}
has 9 linearly independent solutions
\begin{equation}\label{eq 5.6}
{w_1}\left( {x,y} \right) = F_{0;2;2}^{1;0;0}\left[ {\begin{array}{*{20}{c}}
{\,a;}\\
{ - ;}
\end{array}\begin{array}{*{20}{c}}
{\,\,\,\,\,\,\, - ;}\\
{{c_1},{c_2};}
\end{array}\begin{array}{*{20}{c}}
{\,\,\,\,\,\,\,\, - ;}\\
{{d_1},{d_2};}
\end{array}x,y} \right],
\end{equation}
\begin{equation}\label{eq 5.7}
{w_2}\left( {x,y} \right) = {y^{1 - {d_1}}}F_{0;2;2}^{1;0;0}\left[ {\begin{array}{*{20}{c}}
{\,1 - {d_1} + a;}\\
{\,\,\,\,\,\,\,\,\,\,\,\,\,\,\,\,\,\,\,\,\,\, - ;}
\end{array}\begin{array}{*{20}{c}}
{\,\,\,\,\,\,\,\, - ;}\\
{{c_1},{c_2};}
\end{array}\begin{array}{*{20}{c}}
{\,\,\,\,\,\,\,\,\,\,\,\,\,\,\,\,\,\,\,\,\,\,\,\,\,\,\,\,\,\,\,\,\,\,\,\,\,\,\,\,\,\,\,\, - ;}\\
{2 - {d_1},{d_2} - {d_1} + 1;}
\end{array}x,y} \right],
\end{equation}
\begin{equation}\label{eq 5.8}
{w_3}\left( {x,y} \right) = {y^{1 - {d_2}}}F_{0;2;2}^{1;0;0}\left[ {\begin{array}{*{20}{c}}
{\,1 - {d_2} + a;}\\
{\,\,\,\,\,\,\,\,\,\,\,\,\,\,\,\,\,\,\,\,\,\, - ;}
\end{array}\begin{array}{*{20}{c}}
{\,\,\,\,\,\,\,\, - ;}\\
{{c_1},{c_2};}
\end{array}\begin{array}{*{20}{c}}
{\,\,\,\,\,\,\,\,\,\,\,\,\,\,\,\,\,\,\,\,\,\,\,\,\,\,\,\,\,\,\,\,\,\,\,\,\,\,\,\,\,\, - ;}\\
{2 - {d_2},{d_1} - {d_2} + 1;}
\end{array}x,y} \right],
\end{equation}
\begin{equation}\label{eq 5.9}
{w_4}\left( {x,y} \right) = {x^{1 - {c_1}}}F_{0;2;2}^{1;0;0}\left[ {\begin{array}{*{20}{c}}
{\,1 - {c_1} + a;}\\
{\,\,\,\,\,\,\,\,\,\,\,\,\,\,\,\,\,\,\,\,\,\, - ;}
\end{array}\begin{array}{*{20}{c}}
{\,\,\,\,\,\,\,\,\,\,\,\,\,\,\,\,\,\,\,\,\,\,\,\,\,\,\,\,\,\,\,\,\,\,\,\,\,\,\,\,\, - ;}\\
{2 - {c_1},{c_2} - {c_1} + 1;}
\end{array}\begin{array}{*{20}{c}}
{\,\,\,\,\,\,\,\,\, - ;}\\
{{d_1},{d_2};}
\end{array}x,y} \right],
\end{equation}
\begin{equation}\label{eq 5.10}
\begin{array}{l}
{w_5}\left( {x,y} \right) = {x^{1 - {c_1}}}{y^{1 - {d_1}}}\\
 \times F_{0;2;2}^{1;0;0}\left[ {\begin{array}{*{20}{c}}
{\,{c_1} + {d_1} - 2 - a;}\\
{\,\,\,\,\,\,\,\,\,\,\,\,\,\,\,\,\,\,\,\,\,\,\,\,\,\,\,\,\,\,\,\,\,\,\,\, - ;}
\end{array}\begin{array}{*{20}{c}}
{\,\,\,\,\,\,\,\,\,\,\,\,\,\,\,\,\,\,\,\,\,\,\,\,\,\,\,\,\,\,\,\,\,\,\,\,\,\,\,\,\,\, - ;}\\
{2 - {c_1},1 + {c_2} - {c_1};}
\end{array}\begin{array}{*{20}{c}}
{\,\,\,\,\,\,\,\,\,\,\,\,\,\,\,\,\,\,\,\,\,\,\,\,\,\,\,\,\,\,\,\,\,\,\,\,\,\,\,\,\,\,\, - ;}\\
{2 - {d_1},1 + {d_2} - {d_1};}
\end{array}x,y} \right]
\end{array},
\end{equation}
\begin{equation}\label{eq 5.11}
\begin{array}{l}
{w_6}\left( {x,y} \right) = {x^{1 - {c_1}}}{y^{1 - {d_2}}}\\
 \times F_{0;2;2}^{1;0;0}\left[ {\begin{array}{*{20}{c}}
{\,{c_1} + {d_2} - 2 - a;}\\
{\,\,\,\,\,\,\,\,\,\,\,\,\,\,\,\,\,\,\,\,\,\,\,\,\,\,\,\,\,\,\,\,\,\,\, - ;}
\end{array}\begin{array}{*{20}{c}}
{\,\,\,\,\,\,\,\,\,\,\,\,\,\,\,\,\,\,\,\,\,\,\,\,\,\,\,\,\,\,\,\,\,\,\,\,\,\,\,\,\, - ;}\\
{2 - {c_1},1 + {c_2} - {c_1};}
\end{array}\begin{array}{*{20}{c}}
{\,\,\,\,\,\,\,\,\,\,\,\,\,\,\,\,\,\,\,\,\,\,\,\,\,\,\,\,\,\,\,\,\,\,\,\,\,\,\,\,\,\,\, - ;}\\
{1 + {d_1} - {d_2},2 - {d_2};}
\end{array}x,y} \right]
\end{array},
\end{equation}
\begin{equation}\label{eq 5.12}
{w_7}\left( {x,y} \right) = {x^{1 - {c_2}}}F_{0;2;2}^{1;0;0}\left[ {\begin{array}{*{20}{c}}
{\,1 - {c_2} + a;}\\
{\,\,\,\,\,\,\,\,\,\,\,\,\,\,\,\,\,\,\,\,\,\, - ;}
\end{array}\begin{array}{*{20}{c}}
{\,\,\,\,\,\,\,\,\,\,\,\,\,\,\,\,\,\,\,\,\,\,\,\,\,\,\,\,\,\,\,\,\,\,\,\,\,\,\,\, - ;}\\
{{c_1} - {c_2} + 1,2 - {c_2};}
\end{array}\begin{array}{*{20}{c}}
{\,\,\,\,\,\,\,\, - ;}\\
{{d_1},{d_2};}
\end{array}x,y} \right],
\end{equation}
\begin{equation}\label{eq 5.13}
\begin{array}{l}
{w_8}\left( {x,y} \right) = {x^{1 - {c_2}}}{y^{1 - {d_1}}}\\
 \times F_{0;2;2}^{1;0;0}\left[ {\begin{array}{*{20}{c}}
{\,{c_2} + {d_1} - 2 - a;}\\
{\,\,\,\,\,\,\,\,\,\,\,\,\,\,\,\,\,\,\,\,\,\,\,\,\,\,\,\,\,\,\,\,\,\,\, - ;}
\end{array}\begin{array}{*{20}{c}}
{\,\,\,\,\,\,\,\,\,\,\,\,\,\,\,\,\,\,\,\,\,\,\,\,\,\,\,\,\,\,\,\,\,\,\,\,\,\,\,\,\, - ;}\\
{1 + {c_1} - {c_2},2 - {c_2};}
\end{array}\begin{array}{*{20}{c}}
{\,\,\,\,\,\,\,\,\,\,\,\,\,\,\,\,\,\,\,\,\,\,\,\,\,\,\,\,\,\,\,\,\,\,\,\,\,\,\,\,\,\, - ;}\\
{2 - {d_1},1 + {d_2} - {d_1};}
\end{array}x,y} \right]
\end{array},
\end{equation}
\begin{equation}\label{eq 5.14}
\begin{array}{l}
{w_9}\left( {x,y} \right) = {x^{1 - {c_2}}}{y^{1 - {d_2}}}\\
 \times F_{0;2;2}^{1;0;0}\left[ {\begin{array}{*{20}{c}}
{\,{d_2} + {c_2} - 2 - a;}\\
{\,\,\,\,\,\,\,\,\,\,\,\,\,\,\,\,\,\,\,\,\,\,\,\,\,\,\,\,\,\,\,\,\,\,\, - ;}
\end{array}\begin{array}{*{20}{c}}
{\,\,\,\,\,\,\,\,\,\,\,\,\,\,\,\,\,\,\,\,\,\,\,\,\,\,\,\,\,\,\,\,\,\,\,\,\,\,\,\, - ;}\\
{1 + {c_1} - {c_2},2 - {c_2};}
\end{array}\begin{array}{*{20}{c}}
{\,\,\,\,\,\,\,\,\,\,\,\,\,\,\,\,\,\,\,\,\,\,\,\,\,\,\,\,\,\,\,\,\,\,\,\,\,\,\,\,\,\, - ;}\\
{1 + {d_1} - {d_2},2 - {d_2};}
\end{array}x,y} \right]
\end{array},
\end{equation}
where
\begin{equation}\label{eq 5.15}
F_{l;i;j}^{p;q;k}\left[ {\begin{array}{*{20}{c}}
{\left( {{a_p}} \right);}\\
{\left( {{\alpha _l}} \right);}
\end{array}\begin{array}{*{20}{c}}
{\left( {{b_q}} \right);}\\
{\left( {{\beta _m}} \right);}
\end{array}\begin{array}{*{20}{c}}
{\left( {{c_k}} \right);}\\
{\left( {{\gamma _n}} \right);}
\end{array}x,y} \right] = \sum\limits_{r,s = 0}^\infty  {} \frac{{\prod\limits_{j = 1}^p {} {{\left( {{a_j}} \right)}_{r + s}}\prod\limits_{j = 1}^q {} {{\left( {{b_j}} \right)}_r}\prod\limits_{j = 1}^k {} {{\left( {{c_j}} \right)}_s}}}{{\prod\limits_{j = 1}^l {} {{\left( {{\alpha _j}} \right)}_{r + s}}\prod\limits_{j = 1}^m {} {{\left( {{\beta _j}} \right)}_r}\prod\limits_{j = 1}^n {} {{\left( {{\gamma _j}} \right)}_s}r!s!}}{x^r}{y^s}
\end{equation}
are hypergeometric function of Kampe de Feriet (\cite{Appell:27}). Then, in view of (\ref{eq 5.6}) - (\ref{eq 5.14}), the system of hypergeometric equations (\ref{eq 5.4}) has the following special solutions
$$
{\omega _1}\left( {\xi ,\eta } \right) = F_{0;2;2}^{1;0;0}\left[ {\begin{array}{*{20}{c}}
{\,1;}\\
{ - ;}
\end{array}\begin{array}{*{20}{c}}
{\,\,\,\,\,\,\,\,\,\,\,\,\,\,\,\,\,\,\,\,\,\,\,\,\,\,\,\,\,\,\,\, - ;}\\
\displaystyle {\frac{{2 + \alpha }}{3},\frac{{1 + 2\alpha }}{3};}
\end{array}\begin{array}{*{20}{c}}
{\,\,\,\,\,\,\,\,\,\,\,\,\,\,\,\,\,\,\,\,\,\,\,\,\,\,\,\,\,\,\, - ;}\\
\displaystyle {\frac{{2 + \beta }}{3},\frac{{1 + 2\beta }}{3};}
\end{array}\xi ,\eta } \right],
$$
$$
{\omega _2}\left( {\xi ,\eta } \right) = {\eta ^{\displaystyle \frac{{1 - \beta }}{3}}}F_{0;2;2}^{1;0;0}\left[ {\begin{array}{*{20}{c}}
{\,\displaystyle \frac{{4 - \beta }}{3};}\\
{\,\,\,\,\,\,\,\,\,\, - ;}
\end{array}\begin{array}{*{20}{c}}
{\,\,\,\,\,\,\,\,\,\,\,\,\,\,\,\,\,\,\,\,\,\,\,\,\,\,\,\,\,\,\,\, - ;}\\
\displaystyle {\frac{{2 + \alpha }}{3},\frac{{1 + 2\alpha }}{3};}
\end{array}\begin{array}{*{20}{c}}
{\,\,\,\,\,\,\,\,\,\,\,\,\,\,\,\,\,\,\,\,\,\,\,\,\,\,\,\, - ;}\\
\displaystyle {\frac{{4 - \beta }}{3},\frac{{2 + \beta }}{3};}
\end{array}\xi ,\eta } \right],
$$
$$
{\omega _3}\left( {\xi ,\eta } \right) = {\eta ^{\displaystyle \frac{{2\left( {1 - \beta } \right)}}{3}}}F_{0;2;2}^{1;0;0}\left[ {\begin{array}{*{20}{c}}
{\,\displaystyle \frac{{5 - 2\beta }}{3};}\\
{\,\,\,\,\,\,\,\,\,\,\,\,\,\,\, - ;}
\end{array}\begin{array}{*{20}{c}}
{\,\,\,\,\,\,\,\,\,\,\,\,\,\,\,\,\,\,\,\,\,\,\,\,\,\,\,\,\,\,\, - ;}\\
\displaystyle {\frac{{2 + \alpha }}{3},\frac{{1 + 2\alpha }}{3};}
\end{array}\begin{array}{*{20}{c}}
{\,\,\,\,\,\,\,\,\,\,\,\,\,\,\,\,\,\,\,\,\,\,\,\,\,\,\,\,\,\,\, - ;}\\
\displaystyle {\frac{{5 - 2\beta }}{3},\frac{{4 - \beta }}{3};}
\end{array}\xi ,\eta } \right],
$$
$$
{\omega _4}\left( {\xi ,\eta } \right) = {\xi ^{\displaystyle \frac{{1 - \alpha }}{3}}}F_{0;2;2}^{1;0;0}\left[ {\begin{array}{*{20}{c}}
{\,\frac{{4 - \alpha }}{3};}\\
{\,\,\,\,\,\,\, - ;}
\end{array}\begin{array}{*{20}{c}}
{\,\,\,\,\,\,\,\,\,\,\,\,\,\,\,\,\,\,\,\,\,\,\,\,\,\,\,\, - ;}\\
\displaystyle {\frac{{4 - \alpha }}{3},\frac{{2 + \alpha }}{3};}
\end{array}\begin{array}{*{20}{c}}
{\,\,\,\,\,\,\,\,\,\,\,\,\,\,\,\,\,\,\,\,\,\,\,\,\,\,\,\,\,\,\, - ;}\\
\displaystyle {\frac{{2 + \beta }}{3},\frac{{1 + 2\beta }}{3};}
\end{array}\xi ,\eta } \right],
$$
$$
{\omega _5}\left( {\xi ,\eta } \right) = {\xi ^{\displaystyle  \frac{{1 - \alpha }}{3}}}{\eta ^{\displaystyle  \frac{{1 - \beta }}{3}}}F_{0;2;2}^{1;0;0}\left[ {\begin{array}{*{20}{c}}
{\,\displaystyle \frac{{\alpha  + \beta  - 5}}{3};}\\
{\,\,\,\,\,\,\,\,\,\,\,\,\,\,\,\,\,\,\,\,\,\,\, - ;}
\end{array}\begin{array}{*{20}{c}}
{\,\,\,\,\,\,\,\,\,\,\,\,\,\,\,\,\,\,\,\,\,\,\,\,\,\,\,\,\, - ;}\\
\displaystyle  {\frac{{4 - \alpha }}{3},\frac{{2 + \alpha }}{3};}
\end{array}\begin{array}{*{20}{c}}
{\,\,\,\,\,\,\,\,\,\,\,\,\,\,\,\,\,\,\,\,\,\,\,\,\,\,\,\, - ;}\\
\displaystyle  {\frac{{4 - \beta }}{3},\frac{{2 + \beta }}{3};}
\end{array}\xi ,\eta } \right],
$$
$$
{\omega _6}\left( {\xi ,\eta } \right) = {\xi ^{\displaystyle \frac{{1 - \alpha }}{3}}}{\eta ^{\displaystyle \frac{{2\left( {1 - \beta } \right)}}{3}}}F_{0;2;2}^{1;0;0}\left[ {\begin{array}{*{20}{c}}
{\,\displaystyle \frac{{\alpha  + 2\beta  - 6}}{3};}\\
{\,\,\,\,\,\,\,\,\,\,\,\,\,\,\,\,\,\,\,\,\,\,\,\,\,\, - ;}
\end{array}\begin{array}{*{20}{c}}
{\,\,\,\,\,\,\,\,\,\,\,\,\,\,\,\,\,\,\,\,\,\,\,\,\,\,\,\,\, - ;}\\
\displaystyle {\frac{{4 - \alpha }}{3},\frac{{2 + \alpha }}{3};}
\end{array}\begin{array}{*{20}{c}}
{\,\,\,\,\,\,\,\,\,\,\,\,\,\,\,\,\,\,\,\,\,\,\,\,\,\,\,\,\,\,\,\, - ;}\\
\displaystyle {\frac{{4 - \beta }}{3},\frac{{5 - 2\beta }}{3};}
\end{array}\xi ,\eta } \right],
$$
$$
{\omega _7}\left( {\xi ,\eta } \right) = {\xi ^{\displaystyle \frac{{2\left( {1 - \alpha } \right)}}{3}}}F_{0;2;2}^{1;0;0}\left[ {\begin{array}{*{20}{c}}
{\,\displaystyle \frac{{5 - 2\alpha }}{3};}\\
{\,\,\,\,\,\,\,\,\,\,\,\,\, - ;}
\end{array}\begin{array}{*{20}{c}}
{\,\,\,\,\,\,\,\,\,\,\,\,\,\,\,\,\,\,\,\,\,\,\,\,\,\,\,\,\,\,\,\, - ;}\\
\displaystyle {\frac{{4 - \alpha }}{3},\frac{{5 - 2\alpha }}{3};}
\end{array}\begin{array}{*{20}{c}}
{\,\,\,\,\,\,\,\,\,\,\,\,\,\,\,\,\,\,\,\,\,\,\,\,\,\,\,\,\,\,\,\, - ;}\\
\displaystyle {\frac{{2 + \beta }}{3},\frac{{1 + 2\beta }}{3};}
\end{array}\xi ,\eta } \right],
$$
$$
{\omega _8}\left( {\xi ,\eta } \right) = {\xi ^{\displaystyle \frac{{2\left( {1 - \alpha } \right)}}{3}}}{\eta ^{\displaystyle \frac{{1 - \beta }}{3}}}F_{0;2;2}^{1;0;0}\left[ {\begin{array}{*{20}{c}}
{\,\displaystyle \frac{{2\alpha  + \beta  - 6}}{3};}\\
{\,\,\,\,\,\,\,\,\,\,\,\,\,\,\,\,\,\,\,\,\,\,\,\,\,\, - ;}
\end{array}\begin{array}{*{20}{c}}
{\,\,\,\,\,\,\,\,\,\,\,\,\,\,\,\,\,\,\,\,\,\,\,\,\,\,\,\,\,\,\,\, - ;}\\
\displaystyle {\frac{{4 - \alpha }}{3},\frac{{5 - 2\alpha }}{3};}
\end{array}\begin{array}{*{20}{c}}
{\,\,\,\,\,\,\,\,\,\,\,\,\,\,\,\,\,\,\,\,\,\,\,\,\,\,\,\, - ;}\\
\displaystyle {\frac{{4 - \beta }}{3},\frac{{2 + \beta }}{3};}
\end{array}\xi ,\eta } \right],
$$
$$
{\omega _9}\left( {\xi ,\eta } \right) = {\xi ^{\displaystyle \frac{{2\left( {1 - \alpha } \right)}}{3}}}{\eta ^{\displaystyle \frac{{2\left( {1 - \beta } \right)}}{3}}}F_{0;2;2}^{1;0;0}\left[ {\begin{array}{*{20}{c}}
{\,\displaystyle \frac{{2\alpha  + 2\beta  - 7}}{3};}\\
{\,\,\,\,\,\,\,\,\,\,\,\,\,\,\,\,\,\,\,\,\,\,\,\,\,\,\,\,\, - ;}
\end{array}\begin{array}{*{20}{c}}
{\,\,\,\,\,\,\,\,\,\,\,\,\,\,\,\,\,\,\,\,\,\,\,\,\,\,\,\,\,\,\,\, - ;}\\
\displaystyle {\frac{{4 - \alpha }}{3},\frac{{5 - 2\alpha }}{3};}
\end{array}\begin{array}{*{20}{c}}
{\,\,\,\,\,\,\,\,\,\,\,\,\,\,\,\,\,\,\,\,\,\,\,\,\,\,\,\,\,\,\, - ;}\\
\displaystyle {\frac{{4 - \beta }}{3},\frac{{5 - 2\beta }}{3};}
\end{array}\xi ,\eta } \right].
$$
Multiplying each solution by $P = {\left( {\displaystyle \frac{2}{{k + 1}}{\displaystyle t^{k + 1}}} \right)^{ - 1}}$, we finally get special solutions for equation (\ref{eq 5.1}):
\begin{equation}\label{eq 5.16}
{u_1}\left( {x,y,t} \right) = {\lambda _1}PF_{0;2;2}^{1;0;0}\left[ {\begin{array}{*{20}{c}}
{\,1;}\\
{ - ;}
\end{array}\begin{array}{*{20}{c}}
{\,\,\,\,\,\,\,\,\,\,\,\,\,\,\,\,\,\,\,\,\,\,\,\, - ;}\\
\displaystyle {\frac{{2 + \alpha }}{3},\frac{{1 + 2\alpha }}{3};}
\end{array}\begin{array}{*{20}{c}}
{\,\,\,\,\,\,\,\,\,\,\,\,\,\,\,\,\,\,\,\,\,\,\,\, - ;}\\
\displaystyle {\frac{{2 + \beta }}{3},\frac{{1 + 2\beta }}{3};}
\end{array}\xi ,\eta } \right],
\end{equation}
\begin{equation}\label{eq 5.17}
{u_2}\left( {x,y,t} \right) = {\lambda _2}P{\eta ^{\displaystyle \frac{{1 - \beta }}{3}}}F_{0;2;2}^{1;0;0}\left[ {\begin{array}{*{20}{c}}
{\,\displaystyle \frac{{4 - \beta }}{3};}\\
{\,\,\,\,\,\,\,\,\,\, - ;}
\end{array}\begin{array}{*{20}{c}}
{\,\,\,\,\,\,\,\,\,\,\,\,\,\,\,\,\,\,\,\,\,\,\,\, - ;}\\
\displaystyle {\frac{{2 + \alpha }}{3},\frac{{1 + 2\alpha }}{3};}
\end{array}\begin{array}{*{20}{c}}
{\,\,\,\,\,\,\,\,\,\,\,\,\,\,\,\,\,\,\,\,\,\,\, - ;}\\
\displaystyle {\frac{{4 - \beta }}{3},\frac{{2 + \beta }}{3};}
\end{array}\xi ,\eta } \right],
\end{equation}
\begin{equation}\label{eq 5.18}
{u_3}\left( {x,y,t} \right) = {\lambda _3}P{\eta ^{\displaystyle \frac{{2\left( {1 - \beta } \right)}}{3}}}F_{0;2;2}^{1;0;0}\left[ {\begin{array}{*{20}{c}}
{\,\displaystyle \frac{{5 - 2\beta }}{3};}\\
{\,\,\,\,\,\,\,\,\,\,\,\,\, - ;}
\end{array}\begin{array}{*{20}{c}}
{\,\,\,\,\,\,\,\,\,\,\,\,\,\,\,\,\,\,\,\,\,\,\,\, - ;}\\
\displaystyle {\frac{{2 + \alpha }}{3},\frac{{1 + 2\alpha }}{3};}
\end{array}\begin{array}{*{20}{c}}
{\,\,\,\,\,\,\,\,\,\,\,\,\,\,\,\,\,\,\,\,\,\,\,\,\, - ;}\\
\displaystyle {\frac{{5 - 2\beta }}{3},\frac{{4 - \beta }}{3};}
\end{array}\xi ,\eta } \right],
\end{equation}
\begin{equation}\label{eq 5.19}
{u_4}\left( {x,y,t} \right) = {\lambda _4}P{\xi ^{\displaystyle \frac{{1 - \alpha }}{3}}}F_{0;2;2}^{1;0;0}\left[ {\begin{array}{*{20}{c}}
{\,\displaystyle \frac{{4 - \alpha }}{3};}\\
{\,\,\,\,\,\,\,\,\,\,\, - ;}
\end{array}\begin{array}{*{20}{c}}
{\,\,\,\,\,\,\,\,\,\,\,\,\,\,\,\,\,\,\,\,\,\, - ;}\\
\displaystyle {\frac{{4 - \alpha }}{3},\frac{{2 + \alpha }}{3};}
\end{array}\begin{array}{*{20}{c}}
{\,\,\,\,\,\,\,\,\,\,\,\,\,\,\,\,\,\,\,\,\,\,\,\,\, - ;}\\
\displaystyle {\frac{{2 + \beta }}{3},\frac{{1 + 2\beta }}{3};}
\end{array}\xi ,\eta } \right],
\end{equation}
\begin{equation}\label{eq 5.20}
{u_5}\left( {x,y,t} \right) = {\lambda _5}P{\xi ^{\displaystyle \frac{{1 - \alpha }}{3}}}{\eta ^{\displaystyle \frac{{1 - \beta }}{3}}}F_{0;2;2}^{1;0;0}\left[ {\begin{array}{*{20}{c}}
{\,\displaystyle \frac{{\alpha  + \beta  - 5}}{3};}\\
{\,\,\,\,\,\,\,\,\,\,\,\,\,\,\,\,\,\,\, - ;}
\end{array}\begin{array}{*{20}{c}}
{\,\,\,\,\,\,\,\,\,\,\,\,\,\,\,\,\,\,\,\,\,\,\,\,\, - ;}\\
\displaystyle {\frac{{4 - \alpha }}{3},\frac{{2 + \alpha }}{3};}
\end{array}\begin{array}{*{20}{c}}
{\,\,\,\,\,\,\,\,\,\,\,\,\,\,\,\,\,\,\,\,\,\,\,\,\, - ;}\\
\displaystyle {\frac{{4 - \beta }}{3},\frac{{2 + \beta }}{3};}
\end{array}\xi ,\eta } \right],
\end{equation}
\begin{equation}\label{eq 5.21}
{u_6}\left( {x,y,t} \right) = {\lambda _6}P{\xi ^{\displaystyle \frac{{1 - \alpha }}{3}}}{\eta ^{\displaystyle \frac{{2\left( {1 - \beta } \right)}}{3}}}F_{0;2;2}^{1;0;0}\left[ {\begin{array}{*{20}{c}}
{\,\displaystyle \frac{{\alpha  + 2\beta  - 6}}{3};}\\
{\,\,\,\,\,\,\,\,\,\,\,\,\,\,\,\,\,\,\,\,\,\,\, - ;}
\end{array}\begin{array}{*{20}{c}}
{\,\,\,\,\,\,\,\,\,\,\,\,\,\,\,\,\,\,\,\,\,\,\,\,\, - ;}\\
\displaystyle {\frac{{4 - \alpha }}{3},\frac{{2 + \alpha }}{3};}
\end{array}\begin{array}{*{20}{c}}
{\,\,\,\,\,\,\,\,\,\,\,\,\,\,\,\,\,\,\,\,\,\,\,\,\,\,\,\, - ;}\\
\displaystyle {\frac{{4 - \beta }}{3},\frac{{5 - 2\beta }}{3};}
\end{array}\xi ,\eta } \right],
\end{equation}
\begin{equation}\label{eq 5.22}
{u_7} \left( {x,y,t} \right) = {\lambda _7}P{\xi ^{\displaystyle \frac{{2\left( {1 - \alpha } \right)}}{3}}}F_{0;2;2}^{1;0;0}\left[ {\begin{array}{*{20}{c}}
{\,\displaystyle \frac{{5 - 2\alpha }}{3};}\\
{\,\,\,\,\,\,\,\,\,\,\,\,\, - ;}
\end{array}\begin{array}{*{20}{c}}
{\,\,\,\,\,\,\,\,\,\,\,\,\,\,\,\,\,\,\,\,\,\,\,\,\,\,\, - ;}\\
\displaystyle {\frac{{4 - \alpha }}{3},\frac{{5 - 2\alpha }}{3};}
\end{array}\begin{array}{*{20}{c}}
{\,\,\,\,\,\,\,\,\,\,\,\,\,\,\,\,\,\,\,\,\,\,\,\,\,\,\, - ;}\\
\displaystyle {\frac{{2 + \beta }}{3},\frac{{1 + 2\beta }}{3};}
\end{array}\xi ,\eta } \right],
\end{equation}
\begin{equation}\label{eq 5.23}
{u_8}\left( {x,y,t} \right) = {\lambda _8}P{\xi ^{\displaystyle \frac{{2\left( {1 - \alpha } \right)}}{3}}}{\eta ^{\displaystyle \frac{{1 - \beta }}{3}}}F_{0;2;2}^{1;0;0}\left[ {\begin{array}{*{20}{c}}
{\,\displaystyle \frac{{2\alpha  + \beta  - 6}}{3};}\\
{\,\,\,\,\,\,\,\,\,\,\,\,\,\,\,\,\,\,\,\,\,\,\, - ;}
\end{array}\begin{array}{*{20}{c}}
{\,\,\,\,\,\,\,\,\,\,\,\,\,\,\,\,\,\,\,\,\,\,\,\,\,\,\, - ;}\\
\displaystyle {\frac{{4 - \alpha }}{3},\frac{{5 - 2\alpha }}{3};}
\end{array}\begin{array}{*{20}{c}}
{\,\,\,\,\,\,\,\,\,\,\,\,\,\,\,\,\,\,\,\,\,\,\,\,\, - ;}\\
\displaystyle {\frac{{4 - \beta }}{3},\frac{{2 + \beta }}{3};}
\end{array}\xi ,\eta } \right],
\end{equation}
\begin{equation}\label{eq 5.24}
{u_9}\left( {x,y,t} \right) = {\lambda _9}P{\xi ^{\displaystyle \frac{{2\left( {1 - \alpha } \right)}}{3}}}{\eta ^{\displaystyle \frac{{2\left( {1 - \beta } \right)}}{3}}}F_{0;2;2}^{1;0;0}\left[ {\begin{array}{*{20}{c}}
{\,\displaystyle \frac{{2\alpha  + 2\beta  - 7}}{3};}\\
{\,\,\,\,\,\,\,\,\,\,\,\,\,\,\,\,\,\,\,\,\,\,\,\,\,\, - ;}
\end{array}\begin{array}{*{20}{c}}
{\,\,\,\,\,\,\,\,\,\,\,\,\,\,\,\,\,\,\,\,\,\,\,\,\,\,\, - ;}\\
\displaystyle {\frac{{4 - \alpha }}{3},\frac{{5 - 2\alpha }}{3};}
\end{array}\begin{array}{*{20}{c}}
{\,\,\,\,\,\,\,\,\,\,\,\,\,\,\,\,\,\,\,\,\,\,\,\,\,\,\,\, - ;}\\
\displaystyle {\frac{{4 - \beta }}{3},\frac{{5 - 2\beta }}{3};}
\end{array}\xi ,\eta } \right],
\end{equation}
where $ {{\lambda _1},\ldots,{\lambda _9}}$ are constants.

\section{A fourth-order differential equation with two lines of degeneracy}

In the domain $\Omega  = \left\{ {\left( {x,t} \right):\,\,x > 0,\,\,t > 0} \right\}$ , we consider the equation
\begin{equation}\label{eq 6.1}
Lu \equiv {x^n}{u_t} - {t^k}{u_{xxxx}} = 0,\,\,\,n,k = const > 0.	
\end{equation}
Special solutions of equation (\ref{eq 6.1}) will be sought in the form
\begin{equation}\label{eq 6.2}
u (x,t) = P(t)\omega \left( \sigma  \right),
\end{equation}
where
\begin{equation}\label{eq 6.3}
P = {\left( {\frac{1}{{k + 1}}{t^{k + 1}}} \right)^{ - 1}},\,\,\,\,\,\,\,\,\,\sigma  =  - \frac{{k + 1}}{{{{\left( {n + 4} \right)}^4}{t^{k + 1}}}}{x^{n + 4}}.
\end{equation}
Substituting (\ref{eq 6.2}) into equation (\ref{eq 6.1}), we define
$$
\begin{array}{l}
{t^k}P{\omega _{\sigma \sigma \sigma \sigma }}\sigma _x^4 + 6{t^k}P{\sigma _{xx}}\sigma _x^2{\omega _{\sigma \sigma \sigma }} \\
 + {t^k}\left[ {3P\sigma _{xx}^2 + 4P{\sigma _x}{\sigma _{xxx}}} \right]{\omega _{\sigma \sigma }} + \left[ {{t^k}P{\sigma _{xxxx}} - {x^n}P{\sigma _t}} \right]{\omega _\sigma } - {x^n}{P_t}\omega  = 0.
\end{array}
$$
After some calculations, we have
\begin{equation}\label{eq 6.4}
\begin{array}{l}
{\sigma ^3}{\omega _{\sigma \sigma \sigma \sigma }} + \left( \displaystyle {3 + \frac{{3 + \alpha }}{4} + \frac{{2 + 2\alpha }}{4} + \frac{{1 + 3\alpha }}{4}} \right){\sigma ^2}{\omega _{\sigma \sigma \sigma }} + \displaystyle \left( {1 + \frac{{3 + \alpha }}{4} + \frac{{2 + 2\alpha }}{4} + \frac{{1 + 3\alpha }}{4}  } \right.\\
 +\displaystyle  \left. {\frac{{3 + \alpha }}{4}\frac{{2 + 2\alpha }}{4} + \frac{{3 + \alpha }}{4}\frac{{1 + 3\alpha }}{4} + \frac{{2 + 2\alpha }}{4}\frac{{1 + 3\alpha }}{4}} \right)\sigma {\omega _{\sigma \sigma }} \\
 + \displaystyle \left( {\frac{{3 + \alpha }}{4}\frac{{2 + 2\alpha }}{4}\frac{{1 + 3\alpha }}{4} - \sigma } \right){\omega _\sigma } - \omega  = 0,
\end{array}
\end{equation}
where $\alpha  = n/\left( {n + 4} \right)$.  From the general theory (\cite{Appell:27}) it is known that the equation
\begin{equation}\label{eq 6.5}
\begin{array}{l}
{x^3}{w_{xxxx}} + \left( {3 + {c_1} + {c_2} + {c_3}} \right){x^2}{w_{xxx}} + \\
 + \left( {1 + {c_1} + {c_2} + {c_3} + {c_1}{c_2} + {c_1}{c_3} + {c_2}{c_3}} \right)x{w_{xx}} + \left( {{c_1}{c_2}{c_3} - x} \right){w_x} - aw = 0,
\end{array}
\end{equation}
has the following special solutions
\begin{equation}\label{eq 6.6}
{w_1} = {\lambda _1}{}_1{F_3}\left( {a;{c_1},{c_2},{c_3};x} \right),
\end{equation}
\begin{equation}\label{eq 6.7}
{w_2} = {\lambda _2}{x^{1 - {c_1}}}{}_1{F_3}\left( {1 - {c_1} + a;2 - {c_1},1 + {c_2} - {c_1},1 + {c_3} - {c_1};x} \right),
\end{equation}
\begin{equation}\label{eq 6.8}
{w_3} = {\lambda _3}{x^{1 - {c_2}}}{}_1{F_3}\left( {1 - {c_2} + a;1 + {c_1} - {c_2},2 - {c_2},1 + {c_3} - {c_2};x} \right),
\end{equation}
\begin{equation}\label{eq 6.9}
{w_4} = {\lambda _4}{x^{1 - {c_3}}}{}_1{F_3}\left( {1 - {c_3} + a;1 + {c_1} - {c_3},1 + {c_2} - {c_3},2 - {c_3};x} \right),
\end{equation}
where ${\lambda _i}$ are constants, $i = {1,2,3,4}, $   and
\begin{equation}\label{eq 6.10}
{}_1{F_3}\left( {a;{c_1},{c_2},{c_3};x} \right) = \sum\limits_{m = 0}^\infty  {} \frac{{{{\left( a \right)}_m}}}{{{{\left( {{c_1}} \right)}_m}{{\left( {{c_2}} \right)}_m}{{\left( {{c_3}} \right)}_m}m!}}{x^m}.
\end{equation}
Then from equation (\ref{eq 6.4}), in view of (\ref{eq 6.6}) - (\ref{eq 6.9}), taking into account representation (\ref{eq 6.2}), it is easy to determine special solutions to equation (\ref{eq 6.1}):
\begin{equation}\label{eq 6.11}
{u_1}\left( {x,t} \right) = {\bar \lambda _1}{\left( {\frac{1}{{k + 1}}{t^{k + 1}}} \right)^{ - 1}}{}_1{F_3}\left( {1;\frac{{3 + \alpha }}{4},\frac{{2 + 2\alpha }}{4},\frac{{1 + 3\alpha }}{4}; - \frac{{k + 1}}{{{{\left( {n + 4} \right)}^4}{t^{k + 1}}}}{x^{n + 4}}} \right),
\end{equation}
\begin{equation}\label{eq 6.12}
{u_2}\left( {x,t} \right) = {\bar \lambda _2}{\left( {\frac{1}{{k + 1}}{t^{k + 1}}} \right)^{ - \frac{{9 - \alpha }}{4}}}x{}_0{F_2}\left( {\frac{{3 + \alpha }}{4},\frac{{2 + 2\alpha }}{4}; - \frac{{k + 1}}{{{{\left( {n + 4} \right)}^4}{t^{k + 1}}}}{x^{n + 4}}} \right),
\end{equation}
\begin{equation}\label{eq 6.13}
{u_3}\left( {x,t} \right) = {\bar \lambda _3}{\left( {\frac{1}{{k + 1}}{t^{k + 1}}} \right)^{ - \frac{{6 - 2\alpha }}{4}}}{x^2}{}_0{F_2}\left( {\frac{{5 - \alpha }}{4},\frac{{3 + \alpha }}{4}; - \frac{{k + 1}}{{{{\left( {n + 4} \right)}^4}{t^{k + 1}}}}{x^{n + 4}}} \right),
\end{equation}
\begin{equation}\label{eq 6.14}
{u_4}\left( {x,t} \right) = {\bar \lambda _4}{\left( {\frac{1}{{k + 1}}{t^{k + 1}}} \right)^{ - \frac{{11 - 3\alpha }}{4}}}{x^3}{}_0{F_2}\left( {\frac{{6 - 2\alpha }}{4},\frac{{5 - \alpha }}{4}; - \frac{{k + 1}}{{{{\left( {n + 4} \right)}^4}{t^{k + 1}}}}{x^{n + 4}}} \right),
\end{equation}
where ${\bar \lambda _i}$ are constants, $i = {1,2,3,4} $.

%
%


\begin{thebibliography}{99}

\bibitem{Sedov:1}
\refitem{book}
L.~I. Sedov, \emph{Similarity and dimension methods in mechanics}  (Science, Moscow, 1977). (in Russian)

\bibitem{Landau:2}
\refitem{article}
L.~D. Landau and E.~M. Lifshits,\textquotedblleft Investigation of the flow features using the Euler-Tricomi equation,\textquotedblright\, DAN AN USSR, \textbf{96}, 725--728 (1954),~(in Russian).

\bibitem{Lohofer:3}
\refitem{article}
G. Lohofer, \textquotedblleft Theory of an electromagnetically deviated metal sphere. I: Absorbed power,\textquotedblright\, SIAM J. Appl. Math. \, \textbf{49}, 567-581 (1989).

\bibitem{Niukkanen:4}
 \refitem{article}
 A.~W. Niukkanen,\textquotedblleft Generalised hypergeometric series ${}^NF(x_1, ..., x_N)$ arising in physical and quantum chemical applications,\textquotedblright\, J. Phys. A: Math. Gen. \textbf{16}, 1813--1825 (1983).

\bibitem{Courant:5}
\refitem{book}
R Courant, K. Friedrichs, \emph{Supersonic flow and shock waves}  (New York, 1984).

\bibitem{Bers:6}
\refitem{book}
L. Bers, \emph{Mathematical problems of subsonic and transonic gas dynamics}  (M., IL., 1961). (in Russian).

\bibitem{Frankl:7}
\refitem{book}
 F.~I. Frankl, \emph{Selected Works in Gas Dynamics} (Nauka, Moscow, 1973) (in Russian).

\bibitem{Ergashev:8}
\refitem{article} T.~G. Ergashev, A. Hasanov, \textquotedblleft Fundamental solutions of the bi-axially symmetric Helmholtz equation,\textquotedblright\, Uzbek Mathematical Journal \textbf{1}, 55--64 (2018).

\bibitem{Seilkhanova:9}
\refitem{article} R.~B. Seilkhanova, A. Hasanov, \textquotedblleft Particular solutions of generalized Euler-Poisson-
Darboux equation,\textquotedblright\, Electron. J. Diff. Equ. \textbf{9}, 1--10 (2015).

\bibitem{Hasanov, Rassias,Turaev:10}
\refitem{book}
 A. Hasanov, J. ~M. Rassias and M. Turaev, \emph{Functional Equations, difference Inequalities and ULAM Stability
Notions, Fundamental solution for the generalized Elliptic Gellerstedt Equation} (Nova Science Publishers Inc. NY, USA, 2010), \textbf{6}, 73--83.

\bibitem{Hasanov:11}
\refitem{article}
A. Hasanov and E.~T. Karimov,\textquotedblleft Fundamental solutions for a class of three-dimensional elliptic equations with singular coefficients,\textquotedblright\, Applied Mathematic Letters \textbf{22}, 1828--1832 (2009).

\bibitem{Salakhitdinov:12}
\refitem{article}
 M.~S. Salakhitdinov  and A. Hasanov, \textquotedblleft The Fundamental solution for one class of degenerate elliptic equations,\textquotedblright\, More Progresses in Analysis, Proceedings of the 5-th International ISAAC Congress. World Scientific Publishing Co Pte. Ltd., 521--531 (2009).

\bibitem{Hasanov:13}
\refitem{article} A. Hasanov,\textquotedblleft Fundamental solutions for degenerated elliptic equation with two perpendicular lines of degeneration,\textquotedblright\, International J. of Applied Mathematics and Statistics \textbf{8} 13, 41--49 (2008).

\bibitem{Hasanov:14}
\refitem{article} M.~S. Salakhitdinov, A. Hasanov,\textquotedblleft To the theory of the multidimensional Gellerstedt equation,\textquotedblright\, Uzbek Mathematical Journal \textbf{3}, 95--109 (2007). (in Russian)

\bibitem{Hasanov:15}
\refitem{article}
A. Hasanov,\textquotedblleft Fundamental solutions bi-axially symmetric Helmholtz equation,\textquotedblright\, Complex Variables and Elliptic Equations \textbf{52} (8), 673--683 (2007).

\bibitem{Hasanov:16}
\refitem{article} A. Hasanov,\textquotedblleft Some solutions of generalized Rassiass equation,\textquotedblright\, International J. of Applied Mathematics and Statistics \textbf{8} M07, 20--30 (2007).

\bibitem{Hasanov:17}
\refitem{article} A. Hasanov,\textquotedblleft Fundamental Solutions of two Degenerated Elliptic Equations and
Solutions of Boundary Value Problems in Infinite Area,\textquotedblright\, International J. of Applied Mathematics and Statistics \textbf{8} M07, 87--95 (2007).

\bibitem{Hasanov:18}
\refitem{article} A. Hasanov,\textquotedblleft Fundamental solutions of generalized Helmholtz equation,\textquotedblright\, Reports of Uzbek
Academy of Sciences \textbf{1}, 13--15 (2006).

\bibitem{Ergashev:19}
\refitem{article} T.~G. Ergashev, A. Hasanov,\textquotedblleft Fundamental solutions of the bi-axially symmetric Helmholtz
equation,\textquotedblright\, Uzbek Mathematical Journal \textbf{1}, 55--64 (2018).

\bibitem{Srivastava:20}
\refitem{article} H. ~M. Srivastava, A. Hasanov and Junesang Choi,\textquotedblleft Double-Layer Potentials for a Generalized Bi-Axially Symmetric Helmholtz Equation,\textquotedblright\, Sohag Journal of Mathematics An International Journal \textbf{1}, 2, 1--10 (2015).


\bibitem{Salakhitdinov:21}
\refitem{article} M.~S. Salakhitdinov and A. Hasanov,\textquotedblleft The Dirichlet problem for the generalized bi-axially
symmetric Helmholtz equation,\textquotedblright\, Eurasian Mathematical Journal \textbf{3}, 4, 99--110 (2012).

\bibitem{Salakhitdinov:22}
\refitem{article} M.~S. Salakhitdinov M. S.  and A. Hasanov,\textquotedblleft A solution of the Neumann-Dirichlet boundary value problem for generalized bi-axially symmetric Helmholtz equation,\textquotedblright\, Complex Variables and Elliptic Equations \textbf{4}, 53, 355--364 (2008).

\bibitem{Hasanov:23}
\refitem{article} A. Hasanov,\textquotedblleft The solution of the Cauchy problem for generalized Euler--Poisson--Darboux
equation,\textquotedblright\, International J. of Applied Mathematics and Statistics \textbf{8} M07, 30--44 (2007).

\bibitem{Salakhitdinov:24}
\refitem{article} M. ~S. Salakhitdinov and A. Hasanov,\textquotedblleft The Tricomi problem for an equation of mixed type with a no smooth line of degeneracy,\textquotedblright\, Differentsialnye Uravneniya \textbf{1} 19, 110--119 (1983).

\bibitem{Erdelyi:25}
\refitem{book}
A. Erd\'{e}lyi, W. Magnus, F.~Oberhettinger, and F.~G. Tricomi, \emph{Higher Transcendental Functions.}  \textbf{1},
McGraw-Hill, New York, Toronto and London, (1953).

\bibitem{Karimov:26}
\refitem{article} Sh. Karimov,\textquotedblleft On one method for the solution of an analog of the Cauchy problem for a polycaloric equation with singular Bessel operator,\textquotedblright\, Ukrainian Mathematical Journal \textbf{69} 10, 1593--1606 (1918).

\bibitem{Appell:27}
\refitem{book}
 P. Appell and J.~Kamp\'{e} de F\'{e}riet, \emph{ Fonctions Hypergeometriques et Hyperspheriques; Polynomes
d'Hermite} (Gauthier - Villars, Paris, 1926).

\end{thebibliography}
\end{document}